\documentclass{svmult}

\usepackage{amsmath}
\usepackage{amsfonts,amssymb}
\usepackage{enumerate}
\usepackage{graphicx}

\newcommand{\lvt}{\left|\kern-1.35pt\left|\kern-1.3pt\left|}
\newcommand{\rvt}{\right|\kern-1.3pt\right|\kern-1.35pt\right|}

\newtheorem{thm}{Theorem}

\newtheorem{prop}[thm]{Proposition}

\graphicspath{{./}}

\begin{document}

\title*{Optimal Points for Cubature Rules and Polynomial Interpolation on a Square}
 \titlerunning{Optimal Points for Cubature Rules and Polynomial Interpolation on Square}
\author{Yuan Xu}
\institute{ \at Department of Mathematics, University of Oregon,
    Eugene, Oregon 97403-1222. \\ \email{yuan@uoregon.edu}. \\ The work was supported in part by NSF Grant DMS-1510296.}

\maketitle
\index{Xu, Yuan}

\paragraph{Dedicated to Ian H.~Sloan on the occasion of his 80th birthday.}
\bigskip

\abstract{
The nodes of certain minimal cubature rule are real common zeros of a set of orthogonal polynomials of 
degree $n$. They often consist of a well distributed set of points and interpolation polynomials based on them have 
desired convergence behavior. We report what is known and the theory behind by explaining the situation when the domain 
of integrals is a square. }

\section{Introduction}

A numerical integration rule is a finite linear combination of point evaluations that approximates an integral. The degree  
of precision of such a rule is the highest total degree of polynomials that are evaluated exactly. For a fixed degree of 
precision, the minimal rule uses the smallest number of point evaluations. Finding a minimal rule is a difficult problem 
and the most challenging part lies in identifying the set of nodes used in the rule, which is often a desirable set of 
points for polynomial interpolation. For integration on subsets of the real line, a Gaussian quadrature rule is minimal; its 
nodes are known to be zeros of orthogonal polynomials and polynomial interpolation based on the nodes has 
desired convergence behavior. The problem is far less understood in higher dimension, where we have fewer answers 
and many open questions. The purpose of this paper is to explain the situation when the integral domain is a square on 
the plane, for which we know more than on any other domain.

We can work with any fixed square and will fix our choice as 
$$
\Box := [-1,1]^2
$$ 
throughout the paper. Let $\Pi_n^2$ denote the space of polynomials of (total) degree at most $n$ in two real variables, 
where the total degree means the sum of degrees in both variables. It is known that $\dim \Pi_n^2 = (n+1)(n+2)/2$.  
Let $W$ be a nonnegative weight function on the square. For the integral with respect to $W$, a cubature rule 
of degree of precision $m$ (abbreviated as {\it degree $m$} from now on) is a finite sum, defined below, such that
\begin{equation} \label{eq:cuba}
  \int_{\Box} f(x,y) W(x,y) dxdy  = \sum_{k=1}^N \lambda_k f(x_k,y_k), \qquad \hbox{for all $f \in \Pi_m^2$},
\end{equation}
and there exists at least one $f \in \Pi_{m+1}^2$ such that the equality \eqref{eq:cuba} fails to hold. The integer 
$N$ is the number of nodes. The points $(x_k,y_k) \in \mathbb{R}^2$ are called {\it nodes} and the numbers $\lambda_k$ are called 
{\it weights} of the cubature rule, respectively.  We consider only {\it positive} cubature rules for which $\lambda_k$ are all positive. 

As in the case of one variable, the nodes of a minimal cubature rule are closely related to the zeros of orthogonal polynomials. 
A polynomial $P$ is an orthogonal polynomial of degree $n$ with respect to the weight function $W$ if $P \in \Pi_n^2$ and
$$
  \int_{\Box} P (x,y) Q(x,y) W(x,y) dx dy = 0 \qquad \hbox{for all $Q \in \Pi_{n-1}^2$}.
$$ 
Let $\mathcal{V}_n(W)$ denote the space of orthogonal polynomials of degree $n$. Then 
$$
  \dim \mathcal{V}_n(W) = n+1, 
$$ 
as can be seen by applying the GramÐSchmidt process to $x^n, x^{n-1}y, \ldots, x y^{n-1}, y^n$. However, the structure of 
zeros for polynomials of more than one variable can be complicated and what is needed is the common zeros of a family 
of orthogonal polynomials of degree $n$. A common zero of a family of polynomials is a point which is a zero for every 
polynomial in the family. To be more precise, what we often need is to identify a polynomial ideal, $I$, generated by a 
family of orthogonal polynomials in $\mathcal{V}_n(W)$, so that its variety, $V$, is real and zero-dimensional, and the 
cardinality of $V$ equals 
the codimension of $I$. Given the status of real algebraic geometry, this is difficult in general. Only in a few cases can we 
establish the existence of a minimal, or near minimal, cubature rule and identify its generating polynomial ideal explicitly. The 
nodes of such a cubature rule are good points for polynomial interpolation. Indeed, using the knowledge on orthogonal 
polynomial that vanish on the nodes, it is not difficult to construct a polynomial subspace $\Pi_n^*$, so that the problem of 
finding $p$ such that $p(x_i,y_i) = f(x_i,y_i)$ for all nodes $(x_i,y_i)$ of the cubature rule has a unique solution in $\Pi_n^*$. 
Moreover, this interpolation polynomial is easy to compute and has desirable convergence behavior. The above rough 
description applies to all cubature rules. Restricting to the square allows us to describe the idea and results without 
becoming overly tangled by notations. 

The minimal or near-minimal cubature rules offer highly efficient tools for high-precision computation of integrals. It is unlikely, 
however, that they will become a major tool for numerical integration any time soon, because we do not know how to construct 
them in most cases. Moreover, their usage is likely restricted to lower dimension integrals, since they are even less understood
in higher dimensions, where the difficulty increases rapidly as the dimension goes up, and, one could also add, truly 
high-dimensional 
numerical integration is really a different problem (see, for example, \cite{DFS}). Nevertheless, with their deep connection  
to other fields in mathematics and their promise as high dimensional substitute for Gaussian quadrature rules, minimal cubature 
rules are a fascinating object to study. It is our hope that this paper will help attract researchers into this topic.

The paper is organized as follows. We review the theoretic results in the following section. In Section 3, we discuss 
minimal and near minimal cubature rules for the Chebyshev weight functions on the square, which includes a discussion
on the Padua points. In Section 4, we discuss more recent extensions of the results in previous section to a family of weight 
functions that have a singularity on the diagonal of the square. Finally, in Section 5, we describe how cubature rules of lower 
degrees can be established for unit weight function on the square. 

\section{Cubature Rules and Interpolation}

We are interested in integrals with respect to a fixed weight function $W$ over the square, as in \eqref{eq:cuba}, and we assume 
that all moments of $W$ are finite. A typical example of $W$ is the product weight function
$$
  W_{\alpha,\beta}(x,y) := (1-x^2)^\alpha (1-y^2)^\beta, \qquad \alpha,\beta  > -1.
$$ 
This weight function is centrally symmetric, which means that it is symmetric with respect to the origin; more precisely, it 
satisfies $W(x,y) = W(-x,-y)$. If we replace $(1-x^2)^{\alpha}$ by $(1-x)^\alpha(1+x)^\gamma$, with $\gamma \ne \alpha$, 
the resulting weight function will not be centrally symmetric. 

Many of the results below hold for cubature rules with respect to integrals on all domains in the plane, not just 
for the square. We start with the first lower bound for the number of nodes of cubature rules \cite{St}. 

\begin{thm}
Let $n$ be a positive integer and let $m = 2n-1$ or $2n-2$. If the cubature rule \eqref{eq:cuba} is of degree $m$, then 
its number of nodes satisfies 
\begin{equation} \label{eq:1st-lwbd}
        N \ge \dim \Pi_{n-1}^2 = \frac{n(n+1)}{2}. 
\end{equation}
\end{thm}

A cubature rule of degree $m$ is called Gaussian if the lower bound \eqref{eq:1st-lwbd} is attained. In the one-dimensional
case, it is well--known that the Gaussian quadrature rule of degree $2n-1$ has $n = \dim \Pi_{n-1}$ nodes, where 
$\Pi_n$ denote the space of polynomials of degree at most $n$ in one variable, and the same number of nodes is 
needed for the quadrature rule of degree $2n-2$. 

For $n =0,1,2,\ldots$, let $\{P_k^n: 0 \le k \le n\}$ be a basis of $\mathcal{V}_n(W)$. We denote by $\mathbb{P}_n$ the set of this basis 
and we also regard $\mathbb{P}_n$ as a column vector 
$$
  \mathbb{P}_n = (P_0^n, P_1^n, \ldots, P_n^n)^\mathsf{t}, 
$$
where the superscript $\mathsf{t}$ denotes the transpose. The Gaussian cubature rules can be characterized as follows: 

\begin{thm} \label{thm:Gaussian}
Let $\mathbb{P}_s$ be a basis of $\mathcal{V}_s(W)$ for $s = n$ and $n-1$. Then 
\begin{enumerate}
\item A Gaussian cubature rule \eqref{eq:cuba} of degree $2n-1$ exists if, and only if, its nodes are common zeros of 
the polynomials in $\mathbb{P}_n$;
\item A Gaussian cubature rule \eqref{eq:cuba} of degree $2n-2$ exists if, and only if, its nodes are common zeros of
the polynomials in 
$$
       \mathbb{P}_n + \Gamma\, \mathbb{P}_{n-1}, 
$$
where $\Gamma$ is a real matrix of size $(n+1) \times n$. 
\end{enumerate}
\end{thm}
 
For $m=2n-1$, the characterization is classical and established in \cite{My1}; see also \cite{DX,My}. For $m=2n-2$, the 
characterization was established in \cite{MP, Sch}. As in the classical Gaussian quadrature rules, a Gaussian
cubature rule, if it exists, can be derived from integrating the Lagrange interpolation based on its nodes. 

Let $(x_k,y_k): 1 \le k \le \dim \Pi_{n-1}^2$ be distinct points in $\mathbb{R}^2$. The Lagrange interpolation polynomial, 
denoted by $L_n f$, is a polynomial of degree $n$, such that 
$$
      L_n f(x_k,y_k) = f(x_k,y_k), \qquad 1 \le k \le \dim \Pi_{n-1}^2.
$$
If $(x_k,y_k)$ are zeros of a Gaussian cubature rule, then the Lagrange interpolation polynomial is uniquely 
determined. Moreover, let $K_n(\cdot,\cdot)$ be the reproducing kernel of the space $\mathcal{V}_n(W)$, which can be 
written as 
$$
   K_n((x,y),(x',y')) := \sum_{m=0}^n \sum_{k=0}^m P_k^m(x,y) P_k^m(x',y'), 
$$
where $\{P_k^m: 0 \le k \le m\}$ is an orthonormal basis of $\mathcal{V}_m(W)$; then the Lagrange interpolation polynomial 
based on the nodes $(x_k,y_k)$ of the Gaussian cubature rule can be written as 
$$
   L_n f(x,y) = \sum_{k=0}^N f(x_k,y_k) \ell_{k,n}(x,y), \qquad \ell_{k,n} : = \frac{K_{n-1}( (x,y), (x_k,y_k))}{K_{n-1}( (x_k,y_k), (x_k,y_k))}, 
$$
where $\lambda_{k}$ are the cubature weights; moreover, $\lambda_{k} = 1/K_{n-1}( (x_k,y_k), (x_k,y_k))$ is 
clearly positive. 

Another characterization, more explicit, of the Gaussian cubature rules can be given in terms of the coefficient
matrices of the three-term relations satisfied by the orthogonal polynomials. 

For $n =0,1,2,\ldots$, let $\{P_k^n: 0 \le k \le n\}$ be an orthonormal basis of $\mathcal{V}_n(W)$. Then there exist matrices 
$A_{n,i}: (n+1)\times (n+2) $ and $B_{n,i}: (n+1) \times (n+1)$ such that (\cite{DX}),
\begin{equation} \label{eq:3-term}
  x_i \mathbb{P}_n(x) = A_{n,i} \mathbb{P}_{n+1}(x) + B_{n,i}\mathbb{P}_n(x) + A_{n-1,i}^{\mathsf{t}} \mathbb{P}_{n-1}(x), \quad x = (x_1,x_2), 
\end{equation} 
for $i =1,2$. The coefficient matrices $B_{n,i}$ are necessarily symmetric. Furthermore, it is known that $B_{n,i} = 0$ if $W$ is 
centrally symmetric. 

\begin{thm}
Let  $n \in \mathbb{N}$, Assume that the cubature rule \eqref{eq:cuba} is of degree $2n-1$. 
\begin{enumerate}
\item The number of nodes of the cubature rule satisfies 
\begin{equation} \label{eq:2ed-lwbd}
        N \ge \dim \Pi_{n-1}^2 + \frac12 \mathrm{rank} (A_{n-1,1} A_{n-1,2}^\mathsf{t} - A_{n-1,2} A_{n-1,1}^\mathsf{t}).
\end{equation}
\item The cubature is Gaussian if, and only if, $A_{n-1,1} A_{n-1,2}^\mathsf{t} =A_{n-1,2} A_{n-1,1}^\mathsf{t}$.
\item If $W$ is centrally symmetric, then \eqref{eq:2ed-lwbd} becomes 
\begin{equation} \label{eq:3rd-lwbd}
   N \ge \dim \Pi_{n-1}^2+   \left \lfloor \frac{n}{2} \right\rfloor= \frac{n(n+1)}{2} +\left \lfloor \frac{n}{2} \right\rfloor =: N_{\rm min}. 
\end{equation}
In particular, Gaussian cubature rules do not exist for centrally symmetric weight functions.
\end{enumerate}
\end{thm}

The lower bound \eqref{eq:3rd-lwbd} was established by M\"oller in his thesis (see \cite{M}). The more general lower bound 
\eqref{eq:2ed-lwbd} was established in \cite{X94}, which reduces to \eqref{eq:3rd-lwbd} when $W$ is centrally symmetric. 
The non-existence of the Gaussian cubature rule of degree $2n-1$ for centrally symmetric weight functions motivates the 
consideration of {\it minimal cubature rules}, defined as the cubature rule(s) with the smallest number of nodes among all 
cubature rules of the same degree for the same integral. Evidently, the existence of a minimal cubature rule is a tautology 
of its definition.

Cubature rules of degree $2n-1$ that attain M\"oller's lower bound $ N_{\rm min}$ in \eqref{eq:3rd-lwbd} can be characterize
in terms of common zeros of orthogonal polynomials as well. 

\begin{thm} \label{thm:op-zeros}
Let $W$ be centrally symmetric. A cubature rule of degree $2n-1$ attains M\"oller's lower bound \eqref{eq:3rd-lwbd} if,
and only if, its nodes are common zeros of $(n+1) - \left \lfloor \frac{n}{2} \right \rfloor$ many orthogonal polynomials of 
degree $n$ in $\mathcal{V}_n(W)$. 
\end{thm}

This theorem was established in \cite{M}. In the language of polynomial ideal and variety, we say that the nodes of 
the cubature rule are the variety of a polynomial ideal generated by $\left \lfloor \frac{n+1}{2} \right \rfloor+1$ many 
orthogonal polynomials of degree $n$. More general results of this nature were developed in \cite{X94}, which shows, 
in particular, that a cubature rule of degree $2n-1$ with $N= N_{\rm min} +1$ exists if its nodes are common zeros of  
$\left \lfloor \frac{n+1}{2} \right \rfloor$ many orthogonal polynomials of degree $n$ in $\mathcal{V}_n(W)$. 

These cubature rules can also be derived from integrating their corresponding interpolating polynomials. However, since 
$N_{\rm min}$ is not equal to the dimension of $\Pi_{n-1}^2$, we need to define an appropriate polynomial subspace
in order to guarantee that the Lagrange interpolant is unique. Assume that a cubature rule of degree $2n-1$ with 
$N = N_{\rm min}$ exists. Let $\sigma =\left \lfloor \frac{n}{2} \right \rfloor$ and let $\mathcal{P}_n:=\{P_1, \ldots, P_{n-\sigma}\}$ 
be the set of orthogonal polynomials whose common zeros are the nodes of the cubature rule. We can assume, 
without loss of generality, that these polynomials are mutually orthogonal and they form an orthonormal subset of 
$\mathcal{V}_n(W)$. Let $\mathcal{Q}_n:=\{Q_{1}, \ldots, Q_\sigma\}$ be an orthonormal basis of 
$\mathcal{V}_n(W) \setminus \mathrm{span} \, \mathcal{P}_n$, 
so that $\mathcal{P}_n \cup \mathcal{Q}_n$ is an orthonormal basis of $\mathcal{V}_n(W)$. Then it is shown in \cite{X94} 
that there is a unique polynomial in the space 
\begin{equation}\label{eq:Pi_n^*}
  \Pi_n^*:= \Pi_{n-1}^2 \cup \mathrm{span}\, \mathcal{Q}_n
\end{equation}
that interpolates a generic function $f$ on the nodes of the minimal cubature rule; that is, there is a unique polynomial 
$L_n f \in \Pi_n^*$ such that 
$$
   L_n f(x_k,y_k) = f(x_k,y_k), \qquad 1 \le k \le N_{\rm min}, 
$$
where $(x_k,y_k)$ are zeros of the minimal cubature rule. Furthermore, this polynomial can be written as 
\begin{equation} \label{eq:Lnf}
   L_n f(x,y) = \sum_{k=0}^N f(x_k,y_k) \ell_{k,n}(x,y), \qquad \ell_{k,n} : = \frac{K^*_n( (x,y), (x_k,y_k))}{K^*_n( (x_k,y_k), (x_k,y_k))}, 
\end{equation}
where 
\begin{equation} \label{eq:Kn*}
  K_n^*((x,y),(x',y')) = K_{n-1}((x,y),(x',y')) + \sum_{j=1}^\sigma Q_j(x,y) Q_j(x',y'). 
\end{equation}
Integrating $L_n f$ gives a cubature rule with $N_{\rm min}$ nodes that is exact for all polynomials in $\Pi_{2n-1}^2$ and, in particular, 
$\lambda_{k} = 1/K^*_n( (x_k,y_k), (x_k,y_k))$. Furthermore, the above relation between cubature rules and interpolation 
polynomials hold if $\sigma =  \lfloor \frac{n}{2} \rfloor + 1$ and the cubature rule has $N_{\rm min} +1$ points. 

All our examples are given for cubature rules for centrally symmetric cases. We are interested in cubature rules that either
attain or nearly attain the lower bounds, which means Gaussian cubature of degree $2n-2$ or cubature rules of degree $2n-1$ 
with $N_{\rm min}$ nodes or $N_{\rm min}+1$ nodes. When such a cubature rule exists, the Lagrange interpolation polynomials 
based on its nodes possesses good, close to optimal, approximation behavior.  

Because our main interest lies in the existence of our cubature rules and the convergence behavior of our interpolation 
polynomials, we shall not state cubature weights, $\lambda_k$ in \eqref{eq:cuba}, nor explicit formulas for the interpolation 
polynomials throughout this paper. For all cases that we shall encounter below, these cubature weights can be stated 
explicitly in terms of known quantities and interpolation polynomials can be written down in closed forms, which can be found in the references 
that we provide. 

\section{Results for Chebyshev Weight Function}

We start with the product Gegenbauer weight function defined on $[-1,1]^2$ by 
$$
  W_\lambda(x,y) = (1-x^2)^{\lambda-1/2} (1-y^2)^{\lambda-1/2}, \qquad \lambda > -1/2. 
$$ 
The cases $\lambda =0$ and $\lambda =1$ are the Chebyshev weight functions of the first and the second kind, respectively. One mutually
orthogonal basis of $\mathcal{V}_n(W)$ is given by 
$$
   P_k^n(x,y) := C_{n-k}^\lambda(x) C_{k}^\lambda(y), \qquad 0 \le k \le n, 
$$
where $C_n^\lambda$ denotes the usual Gegenbauer polynomial of degree $n$. When $\lambda = 0$, $C_n^\lambda$ is replaced by $T_n$, the 
Chebyshev polynomial of the first kind, and when $\lambda =1$, $C_n^\lambda = U_n$, the Chebyshev polynomial of the second kind. 
Setting $x = \cos \theta$, we have 
$$
T_n(x) = \cos n \theta \quad \hbox{and} \quad U_n(x) = \frac{\sin (n+1)\theta}{\sin \theta}. 
$$
In the following we always assume that $C_k^\lambda(x) = U_k(x)  = T_k(x)=0$ if $k < 0$. 
 
The first examples of minimal cubature rules were given for Chebyshev weight functions soon after \cite{M}. We start
with the Gaussian cubature rules for Chebyshev weight function of the second type in \cite{MP}.

\begin{thm} 
For the product Chebyshev weight function $W_1$ of the second kind, the Gaussian cubature rules of degree $2n-2$ exist. 
Their nodes can be explicitly given by
\begin{align} \label{eq:nodesU}
\begin{split}
   & (\cos \tfrac{2i \pi}{n+2}, \cos \tfrac{(2j-1) \pi}{n+1}),   \quad  1 \le i \le (n+1)/2, \quad 1 \le j \le (n+1)/2, \\
   & (\cos \tfrac{(2i -1)\pi}{n+2}, \cos \tfrac{2j \pi}{n+1}), \quad   1 \le i \le n/2 + 1, \quad 1 \le j \le n/2,
\end{split}
\end{align}
which are common zeros of the polynomials 
$$
 U_{n-k}(x) U_k(y) - U_k(x) U_{n-1-k}(y), \qquad 0 \le k \le n.
$$
\end{thm} 

However, $W_1$ remains the only weight function on the square for which the Gaussian cubature rules of degree
$2n-2$ are known to exist for all $n$. For other weight functions, for example, the constant weight function
$W_{1/2} (x,y) =1$, the existence is known only for small $n$; see the discussion in the last section. 

For minimal cubature of degree $2n-1$ that attains M\"oller's lower bound, we are in better position. The first result 
is again known for Chebyshev weight functions. 

\begin{thm}
For the product Chebyshev weight function $W_0$ of the first kind, the cubature rules of degree $2n-1$ that attain
the lower bound \eqref{eq:3rd-lwbd} exist. Moreover, for $n=2m$, their nodes can be explicitly given by
\begin{align} \label{eq:nodesT1}
\begin{split}
& (\cos \tfrac{i \pi}{m}, \cos \tfrac{(2j+1)\pi}{2m}), \quad\quad   0 \le i \le m, \quad 0 \le j \le m-1, \\
& (\cos \tfrac{(2 i+1) \pi}{m}, \cos \tfrac{j\pi}{m}),  \quad\quad  0 \le i \le m, \quad 1 \le j \le m, 
\end{split}
\end{align}
which are common zeros of the polynomials 
$$
    T_{2m - k +1}(x)T_{m-1}(y) - T_{k-1}(x) T_{m-k+1}(y), \qquad 1 \le k \le m+1.
$$
\end{thm}

\begin{figure}[htbp]  
\centering
   \includegraphics[width=2.3in]{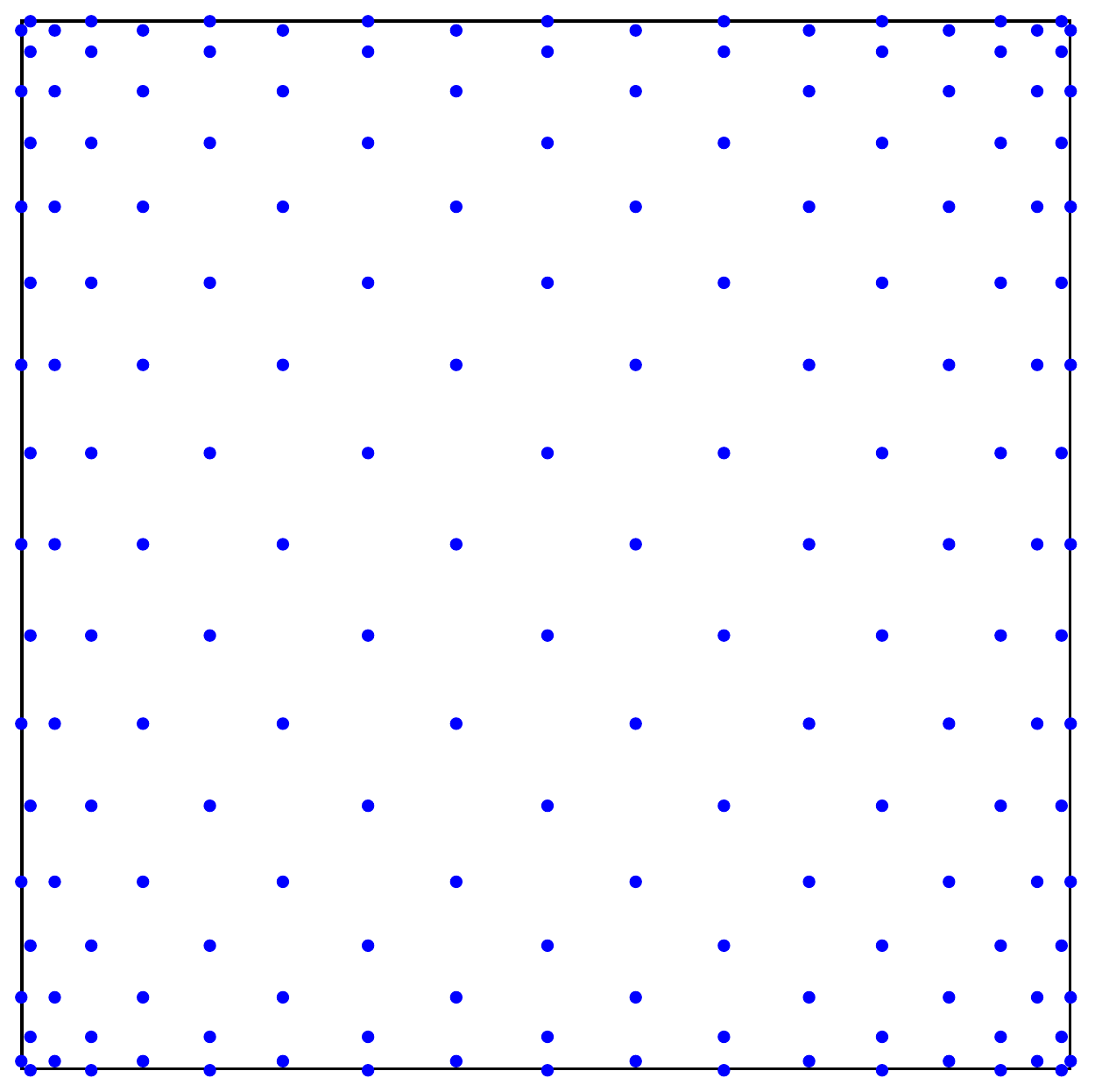} \,\,  \includegraphics[width=2.3in]{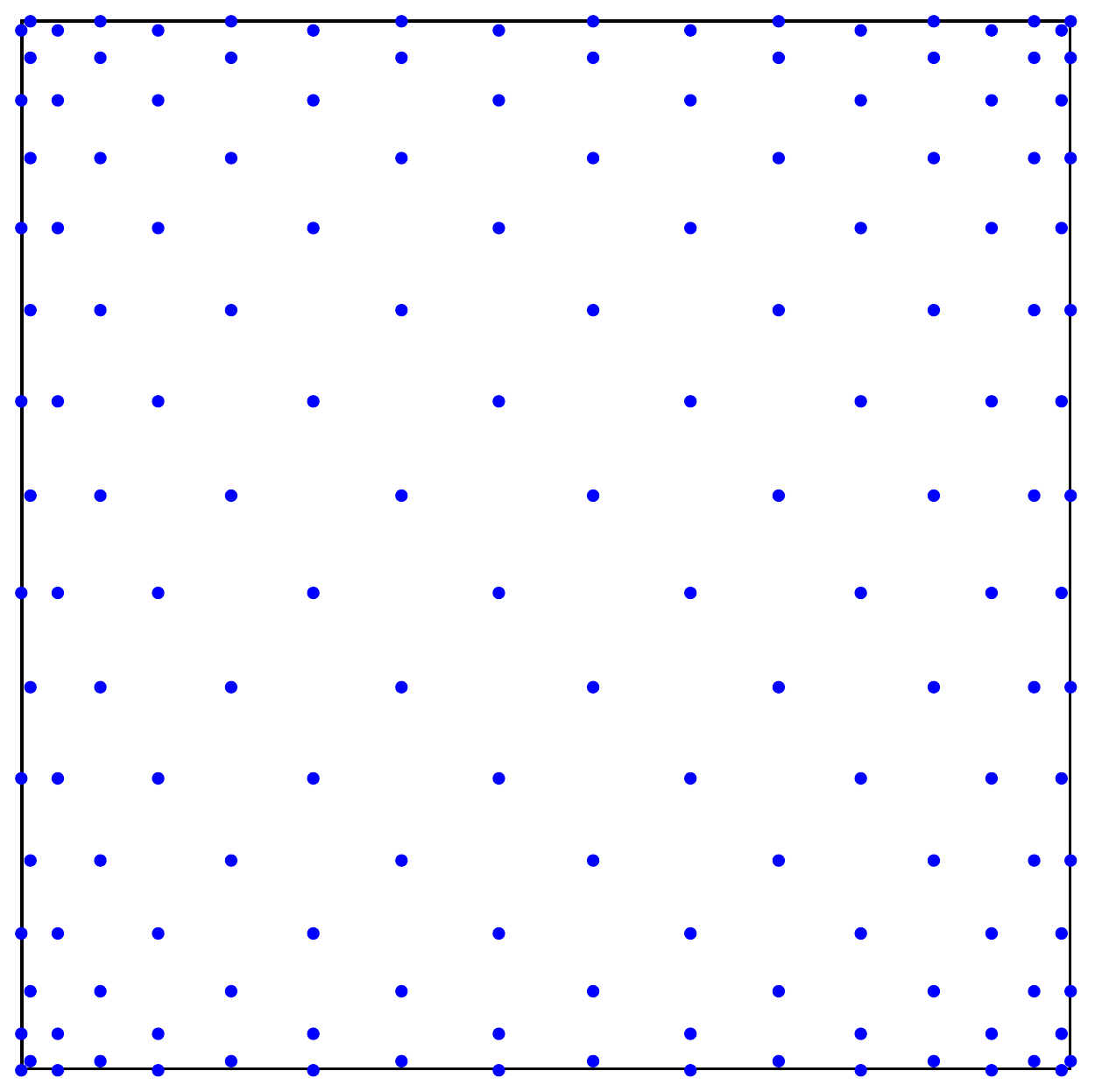} 
   \caption{Left: 180 nodes for minimal cubature rule of degree 35      
    Right: 162 nodes for near-minimal cubature rule of degree 33}
 \label{fig:min}
\end{figure}

For $n=2m$, this was first established in \cite{MP}, using the characterization in \cite{M}, and it was later proved 
by other methods \cite{BP, LSX}. The case for $n=2m-1$ is established more recently in \cite{X16}, for
which the structure of orthogonal polynomials that vanish on the nodes is more complicated, see the discussion 
after Theorem 4.2. The analog of the explicit construction in the case $n=2m$ holds for cubature rules 
of degree $2n-1$, with $n=2m-1$, that have one more node than the lower bound \eqref{eq:3rd-lwbd} \cite{X94}. The 
nodes of these formulas are by \cite{X96}
\begin{align} \label{eq:nodesT2}
\begin{split}
& (\cos \tfrac{2i \pi}{2m-1}, \cos \tfrac{2j\pi}{2m-1}), \quad\quad   0 \le k \le m-1, \, 0 \le j \le m-1, \\
& (\cos \tfrac{(2 m-2i-1) \pi}{2m-1}, \cos \tfrac{(2m-j-1)\pi}{2m-1}),  \quad\quad  0 \le i \le m-1, \, 1 \le j \le m-1, 
\end{split}
\end{align}
and they are common zeros of the polynomials 
$$
    T_{2m - k}(x)T_{k-1}(y) - T_{k-1}(x) T_{2m-k}(y), \qquad 1 \le k \le m.
$$
These points are well distributed. Two examples are depicted in Figure \ref{fig:min}.

The Lagrange interpolation polynomials based on the nodes of these cubature rules were first studied in \cite{X96}. 
Let $L_n f(x,y)$ denote the Lagrange interpolation polynomial based on the nodes \eqref{eq:nodesT1} for 
$n =2m$ and on \eqref{eq:nodesT2} for $n=2m-1$, which belongs to the space $\Pi_n^*$ defined in \eqref{eq:Pi_n^*}. 
Using the Christoffel-Darboux formula in two variables, these interpolation polynomials can be given explicitly. Their 
convergence behavior is about optimal among all interpolation polynomials on the square. To be more precise, we 
introduce the following notation. 

Let $\|\cdot\|_p$ denote the usual $L^p$ norm of the space $L^p(\Box, W_0)$ for $1 \le p < \infty$, and define it
as the uniform norm on the square $\Box$ when $p =\infty$. For $f \in C(\Box)$, let $E_n (f)_\infty$ be the error of 
best approximation by polynomials from $\Pi_n^2$ in the uniform norm; that is, 
$$
  E_n (f)_\infty = \inf_{ P \in \Pi_n^2} \|f - P \|_\infty. 
$$

\begin{thm}\label{thm:conv1}
Let $f$ be a continuous function on $\Box$. Then 
\begin{enumerate}[1.]
\item There is a constant $c > 0$, independent of $n$ and $f$, such that 
$$
\|f - L_n f\|_p \le c \, E_n (f)_\infty, \qquad 1 \le p < \infty; 
$$
\item The Lebesgue constant $\|L_n\|_\infty:= \sup_{\|f\|_\infty \ne 0} \|L_n f\|_\infty$ satisfies 
$$
  \|L_n \|_\infty = \mathcal{O}( (\log n)^2),
$$
which is the optimal order among all projection operators from $C(\Omega) \mapsto \Pi_n^*$. 
\end{enumerate}
\end{thm}

The first item was proved in \cite{X96}, which shows that $L_n f$ behaves like polynomials of best approximation 
in $L^p$ norm when $1 \le p < \infty$. The second one was proved more recently in \cite{BMV}, which gives the upper
bound of the Lebesgue constant; that this upper bound is optimal was established in \cite{SV}. These results indicate
that the set of points \eqref{eq:nodesT1} is optimal for both numerical integration and interpolation. These interpolation 
polynomials were also considered in \cite{Ha1}, and further extended in \cite{Ha1A, Ha2}, where points for other 
Chebyshev weights \cite{MP}, including $(1-x^2)^{\pm \frac{1}{2}}(1-y^2)^{\mp \frac{1}{2}}$, are considered. 

The interpolation polynomial $L_n f$ defined above is of degree $n$ and its set of interpolation points has the 
cardinality $\dim \Pi_{n-1}^2 + \lfloor n/2\rfloor$ or one more. One could ask if it is possible to identify another set of 
points, say $X_n$, that has the cardinality $\dim \Pi_n^2$ and is just as good, which means that the interpolation 
polynomials based on $X_n$ should have the same convergence behavior as those in Theorem \ref{thm:conv1} and 
the cubature rule with $X_n$ as the set of nodes should be of the degree of precision $2n-1$. If such an $X_n$ exists, 
the points in $X_n$ need to be common zeros of polynomials of the form 
$$
\mathbb{P}_{n+1} + \Gamma_1 \mathbb{P}_n + \Gamma_2 \mathbb{P}_{n-1},
$$ 
where $\Gamma_1$ and $\Gamma_2$ are matrices of sizes $(n+2) \times (n+1)$ and $(n+2) \times n$, respectively. 
For $W_0$, such a set indeed exists and known as the Padua points \cite{Pad1, Pad2}. One version of these points is  
\begin{align} \label{eq:nodesPad}
\begin{split}
  X_n: =  & \left \{ (\cos \tfrac{2i \pi}{n}, \cos \tfrac{(2j-1) \pi}{n+1}),   \quad  0 \le i \le \lfloor \tfrac{n}{2} \rfloor \quad 1 \le j \le  \lfloor \tfrac{n}{2} \rfloor +1,  \right.  \\
   & \left. (\cos \tfrac{(2i -1)\pi}{n}, \cos \tfrac{ (2j -2)\pi}{n+1}), \quad   1 \le i \le \lfloor \tfrac{n}{2} \rfloor+1,
      \quad 1 \le  j \le \lfloor \tfrac{n}{2} \rfloor +2\right \},
\end{split}
\end{align}
which are common zeros of polynomials $Q_k^{n+1}$, $0 \le k \le n+1$, defined by 
\begin{align} \label{Q0}
Q_0^{n+1}(x, y)& =T_{n+1}(x)-T_{n-1}(x), \\
  Q_k^{n+1}(x,y)& =T_{n-k+1}(x)T_k(y)+T_{n-k+1}(y)T_{k-1}(x), \quad 1 \le k \le n+1. 
\end{align}

\begin{thm}
For $n \in \mathbb{N}$, let $X_n$ be defined as in \eqref{eq:nodesPad}. Then $|X_n| = \dim \Pi_n^2$ and
\begin{enumerate}[    1.]
\item There is a cubature rule of degree $2n-1$ with $X_n$ as its set of nodes.
\item There is a unique polynomial of degree $n$ that interpolates at the points in $X_n$, which enjoys the 
same convergence as that of $L_n f$ given in Theorem \ref{thm:conv1}.
\end{enumerate}
\end{thm}

\begin{figure}[htbp]  
\centering
 \includegraphics[width=2.5in]{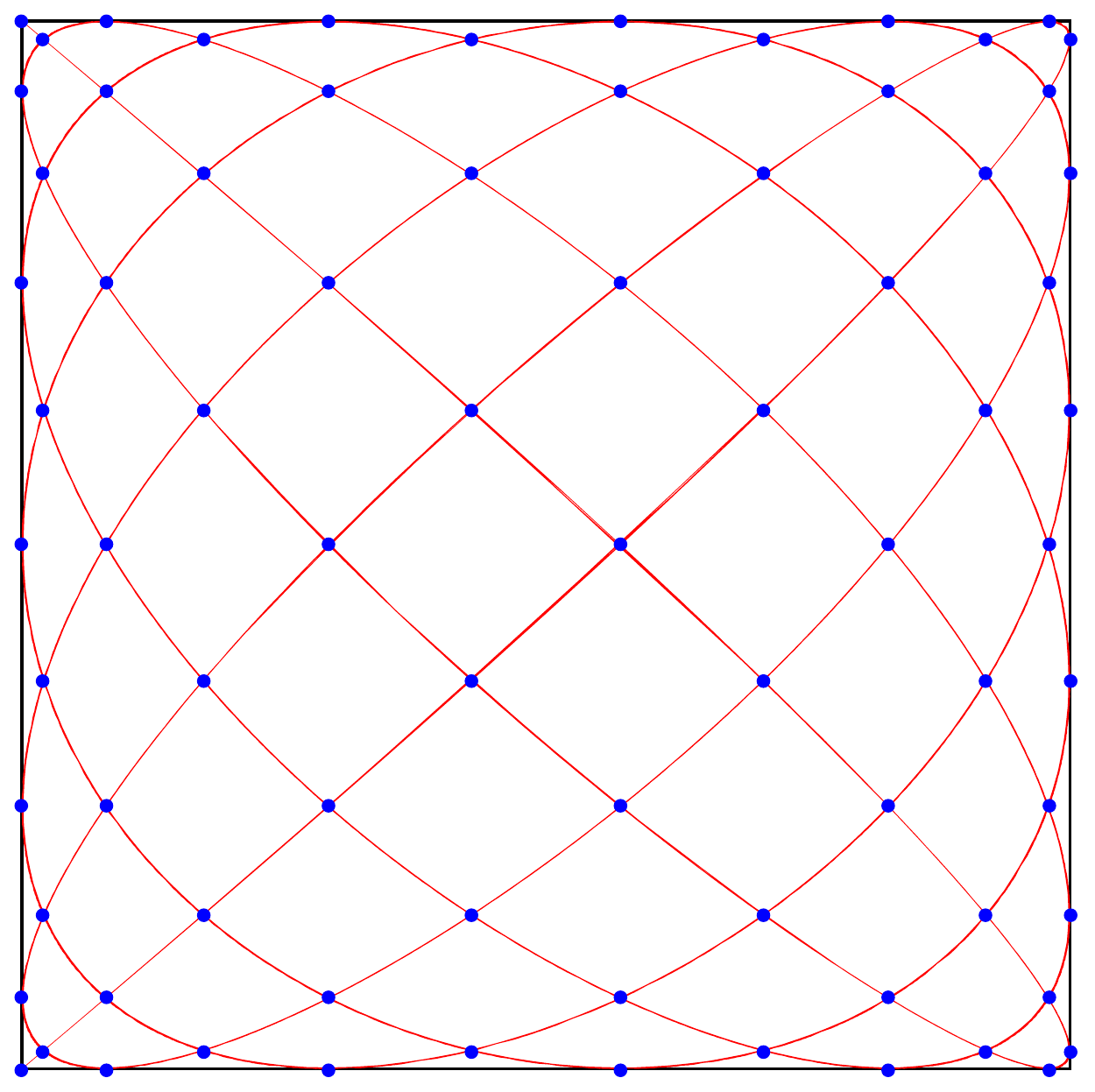}  
   \caption{$78$ Padua points ($n=11$) and their generating curve}
 \label{fig:Padua}
\end{figure}

One interesting property of the Padua points is that they are self-intersection points of a Lissajous curve. For $X_n$ 
given in \eqref{eq:nodesPad}, the curve is given by $Q_0^{n+1}$ or, in parametric form,
$$
     (- \cos ((n + 1) t),  - \cos (n t)), \qquad  0 \le t \le (2 n + 1) \pi, 
$$
as shown in Figure \ref{fig:Padua}. The generating curve offers a convenient tool for studying the interpolation 
polynomial based on Padua points.  
 
More generally, a Lissajous curve takes of the form $(\cos ((n+p)t), \cos (nt))$ with positive integers $n$ and $p$ such
that $n$ and $n+p$ are relatively prime. It is known \cite{F} that such a curve has $(n-1)(n+p-1)/2$ self-intersection points  
inside $[-1,1]^2$. For $p \ne 1$, the number is not equal to the full dimension of $\Pi_m$ for any $m$ in general. 
Nonetheless, these points turn out to be good points for cubature rules and for polynomial interpolation, as shown in \cite{E1,E2}. 

\section{Results for a Family of Weight Functions}

In this section we consider a family of weight functions that include the Chebyshev weight functions as special cases. 
Let $w$ be a weight function on the interval $[-1,1]$. For $\gamma > -1/2$, we define a weight function
\begin{align*}
    \mathcal{W}_\gamma (x,y)  := \, &  w(\cos (\theta-\phi)) w(\cos  (\theta+\phi)) |x^2 - y^2| (1-x^2)^\gamma (1-y^2)^\gamma,   \\ 
         & \qquad \hbox{where} \, \, x = \cos \theta, \,\, y = \cos \phi, \quad
           (x,y) \in [-1,1]^2. \notag
\end{align*}
When $w$ is the Jacobi weight function $w_{\alpha,\beta}(x):= (1-x)^\alpha (1+x)^\beta$, we denote the weight function $\mathcal{W}_\gamma$ by 
$W_{\alpha,\beta,\gamma}$. It is not difficult to verify that 
\begin{equation} \label{eq:Wabg}
   W_{\alpha,\beta, \gamma}(x,y) := |x+y|^{2\alpha +1}  |x- y|^{2\beta +1}  
    (1-x^2)^\gamma (1-y^2)^\gamma.
\end{equation}
In the special cases of $\alpha = \beta = - \frac{1}{2}$ and $\gamma = \pm \frac{1}{2}$,  these are exactly the Chebyshev weight functions. It was 
proved recently in \cite{X12a,X16}, rather surprisingly, that the results in the previous section can be extended to these 
weight functions. First, however, we describe a family of mutually orthogonal polynomials. To be
more precise, we state this basis only for the weight function $W_{\alpha,\beta,\pm \frac{1}{2}}$. 

For $\alpha,\beta > -1$, let $p_n^{(\alpha,\beta)}$ be the normalized Jacobi polynomial of degree $n$, so that 
$c_{\alpha,\beta} \int_{-1}^1 |p_n^{(\alpha,\beta)}(x)|^2 w_{\alpha,\beta}(x) dx =1$ and $p_0^{(\alpha,\beta)}(x) =1$. 
For $x = \cos \theta$ and $y = \cos \phi$, 
we define 
\begin{align*} 
 & P_{k,n}^{\alpha,\beta, -\frac12} (2 x y, x^2+y^2 -1)  \\ 
 &  \, :=   p_n^{(\alpha,\beta)}(\cos (\theta - \phi)) p_k^{(\alpha,\beta)}(\cos (\theta+\phi)) 
           +  p_k^{(\alpha,\beta)} (\cos (\theta - \phi)) p_n^{(\alpha,\beta)}(\cos (\theta+\phi)), \\
 & P_{k,n}^{\alpha,\beta, \frac12} (2 x y, x^2+y^2 -1)  \\ 
 &  \, :=  \frac{p_{n+1}^{(\alpha,\beta)}(\cos (\theta - \phi)) p_k^{(\alpha,\beta)}(\cos (\theta+\phi)) 
           -  p_k^{(\alpha,\beta)} (\cos (\theta - \phi)) p_{n+1}^{(\alpha,\beta)}(\cos (\theta+\phi))}{2 \sin \theta \sin \phi}. 
\end{align*}
It turns out that $P_{k,n}^{\alpha,\beta,\pm \frac{1}{2}}(u,v)$ itself is a polynomial of degree $n$ in the variables $u$ and $v$, 
as can be seen by the elementary trigonometric identities
\begin{equation*} 
    2 x y = \cos (\theta - \phi) + \cos (\theta + \phi) \quad\hbox{and} \quad x^2+y^2 -1 = \cos (\theta -  \phi) \cos (\theta + \phi),  
\end{equation*} 
and the fundamental theorem of symmetric polynomials. Furthermore, $P_{k,n}^{\alpha,\beta,\pm \frac{1}{2}}(u,v)$, first studied 
in \cite{K74a}, are orthogonal polynomials with respect to a weight function on a domain bounded by a parabola and two 
straight lines and the weight function admit Gaussian cubature rules of all degrees \cite{SX}. These polynomials are closely 
related to the orthogonal polynomials in $\mathcal{V}_n(W_{\alpha,\beta, - \frac{1}{2}})$, as shown in the following proposition
established in \cite{X12b}: 

\begin{prop}
Let $\alpha,\beta > -1$. A mutually orthogonal basis for $\mathcal{V}_{2m}(W_{\alpha,\beta, - \frac{1}{2}})$ is given by 
\begin{align} \label{eq:Qeven}
\begin{split}
 {}_1Q_{k,2m}^{\alpha,\beta,\pm \frac{1}{2}}(x,y):= & P_{k,m}^{\alpha,\beta,\pm \frac{1}{2}}(2xy, x^2+y^2 -1), \quad 0 \le k \le m, \\
 {}_2Q_{k,2m}^{\alpha,\beta,\pm \frac{1}{2}}(x,y) := & (x^2-y^2)  P_{k,m-1}^{\alpha+1,\beta+1,\pm \frac{1}{2}}(2xy, x^2+y^2 -1),  
   \quad 0 \le k \le m-1, 
\end{split}
\end{align}
and a mutually orthogonal basis for $\mathcal{V}_{2m+1}(W_{\alpha,\beta,\pm \frac{1}{2}})$ is given by 
\begin{align}  \label{eq:Qodd}
\begin{split}
   {}_1Q_{k,2m+1}^{\alpha,\beta,\pm \frac{1}{2}}(x,y):= &  (x+y)P_{k,m}^{\alpha,\beta+1,\pm \frac{1}{2}}(2xy, x^2+y^2 -1), 
     \quad 0 \le k \le m, \\
   {}_2Q_{k,2m+1}^{\alpha,\beta,\pm \frac{1}{2}}(x,y) := &  (x-y)P_{k,m-1}^{\alpha+1,\beta,\pm \frac{1}{2}}(2xy, x^2+y^2 -1),  
   \quad 0 \le k \le m. 
\end{split}
\end{align}
\end{prop}
 
The orthogonal polynomials in \eqref{eq:Qeven} of degree $2n$ are symmetric polynomials in $x$ and $y$, and they 
are invariant under $(x,y) \mapsto (-x,-y)$. Notice, however, that the product Chebyshev polynomials do not possess 
such symmetries, even though $W_{-\frac{1}{2},-\frac{1}{2}, \pm \frac{1}{2}}$ are the Chebyshev weight functions. 

We now return to cubature rules and interpolation and state the following theorem. 

\begin{thm}
The minimal cubature rules of degree $2n-1$ that attain the lower bound \eqref{eq:3rd-lwbd} exist for the weight function 
$\mathcal{W}_{\pm \frac{1}{2}}$ when $n=2m$. Moreover, the same holds for the weight function $W_{\alpha,\beta,\pm \frac{1}{2}}$ when $n =2m+1$.  
\end{thm}

The orthogonal polynomials whose common zeros are nodes of these cubature rules, as described in Theorem \ref{thm:op-zeros}, 
can be identified explicitly. Let us consider only $W_{\alpha,\beta,-\frac{1}{2}}$. For $n=2m$, these polynomials can be chosen as 
${}_1Q_{k,2m}^{\alpha,\beta, -\frac{1}{2}}$, $0\le k \le m$, in \eqref{eq:Qeven}. For $n = 2m+1$, they can be chosen as 
${}_2Q_{k,2m+1}^{\alpha,\beta, - \frac{1}{2}}$, $0\le k \le m$, in \eqref{eq:Qodd}, together with one more polynomial 
\begin{align*}
  q_m (x,y) = (x+y) & \left[p_m^{(\alpha,\beta+1)}(\cos (\theta-\phi))p_m^{(\alpha+1,\beta)}(\cos (\theta+\phi)) \right. \\
          & \left. + p_m^{(\alpha,\beta+1)}(\cos (\theta+\phi))  p_m^{(\alpha+1,\beta)}(\cos (\theta-\phi))\right] 
\end{align*}
in $\mathcal{V}_{2m+1}(W_{\alpha,\beta,-\frac{1}{2}})$, as shown in \cite{X16}. For $n = 2m$, the nodes of the minimal cubature 
rules for $W_{\alpha,\beta,-\frac{1}{2}}$ are not as explicit as those for $n=2m$. For interpolation, it is often easier to work with the 
near minimal cubature rule of degree $2n-1$ when $n= 2m+1$, whose number of nodes is just one more than the minimal 
number $N_{\min}$ in \eqref{eq:3rd-lwbd}. The nodes of the these near minimal rules are common zeros of 
${}_2Q_{k,2m+1}^{\alpha,\beta,-\frac{1}{2}}$, $0 \le k \le m$, and a quasi-orthogonal polynomial of the form 
${}_1Q_{k,2m+2}^{\alpha,\beta,-\frac{1}{2}} - a_{k,m} {}_1Q_{k,2m}^{\alpha,\beta,-\frac{1}{2}}$, where $a_{k,m}$ are specific constants 
(\cite[Theorem 3.5]{X16}). 

The nodes of these cubature rules can be specified. For $\alpha,\beta > -1$ and $1 \le k \le m$, let $\cos \theta_{k,m}^{\alpha,\beta}$ be 
the zeros of the Jacobi polynomial $P_m^{\alpha,\beta}$ so that 
$$
   0  < \theta_{1,m}^{\alpha,\beta} < \ldots < \theta_{m,m}^{\alpha,\beta} < \pi,
$$ 
and we also define $\theta_{0,m}^{\alpha,\beta} =0$. We further define 
\begin{align*} 
  s_{j,k}^{\alpha,\beta}: =  \cos \tfrac{\theta_{j,n}-\theta_{k,n}}{2} \quad\hbox{and}\quad
                t_{j,k}^{\alpha,\beta}:= \cos \tfrac{\theta_{j,n} + \theta_{k,n}}{2}, \qquad \hbox{where} \quad \theta_{k,n} = \theta_{k,n}^{\alpha,\beta}. 
\end{align*}
For $n =2m$, the nodes of the minimal cubature rule of degree $2n-1$ consist of
\begin{equation*}
  X_{2m}^{\alpha,\beta}:=\{(s_{j,k}^{\alpha,\beta}, t_{j,k}^{\alpha,\beta}), (t_{j,k}^{\alpha,\beta}, s_{j,k}^{\alpha,\beta}), (-s_{j,k}^{\alpha,\beta}, -t_{j,k}^{\alpha,\beta}), 
  (-t_{j,k}^{\alpha,\beta}, -s_{j,k}^{\alpha,\beta}): 1 \le j \le k \le m\}. 
\end{equation*}
For $n=2m+1$, the nodes of the near minimal cubature rule of degree $2n-1$ consist of 
\begin{align*}
  X_{2m+1}^{\alpha,\beta}: =\{(s_{j,k}^{\alpha+1,\beta},  t_{j,k}^{\alpha+1,\beta}), 
     &  (t_{j,k}^{\alpha+1,\beta},  s_{j,k}^{\alpha+1,\beta}),    (-s_{j,k}^{\alpha+1,\beta}, -t_{j,k}^{\alpha+1,\beta}), \\
    & (-t_{j,k}^{\alpha+1,\beta}, -s_{j,k}^{\alpha+1,\beta}):   0 \le j \le k \le m\}. 
\end{align*}

\begin{figure}[htbp]  
\centering
   \includegraphics[width=2.3in]{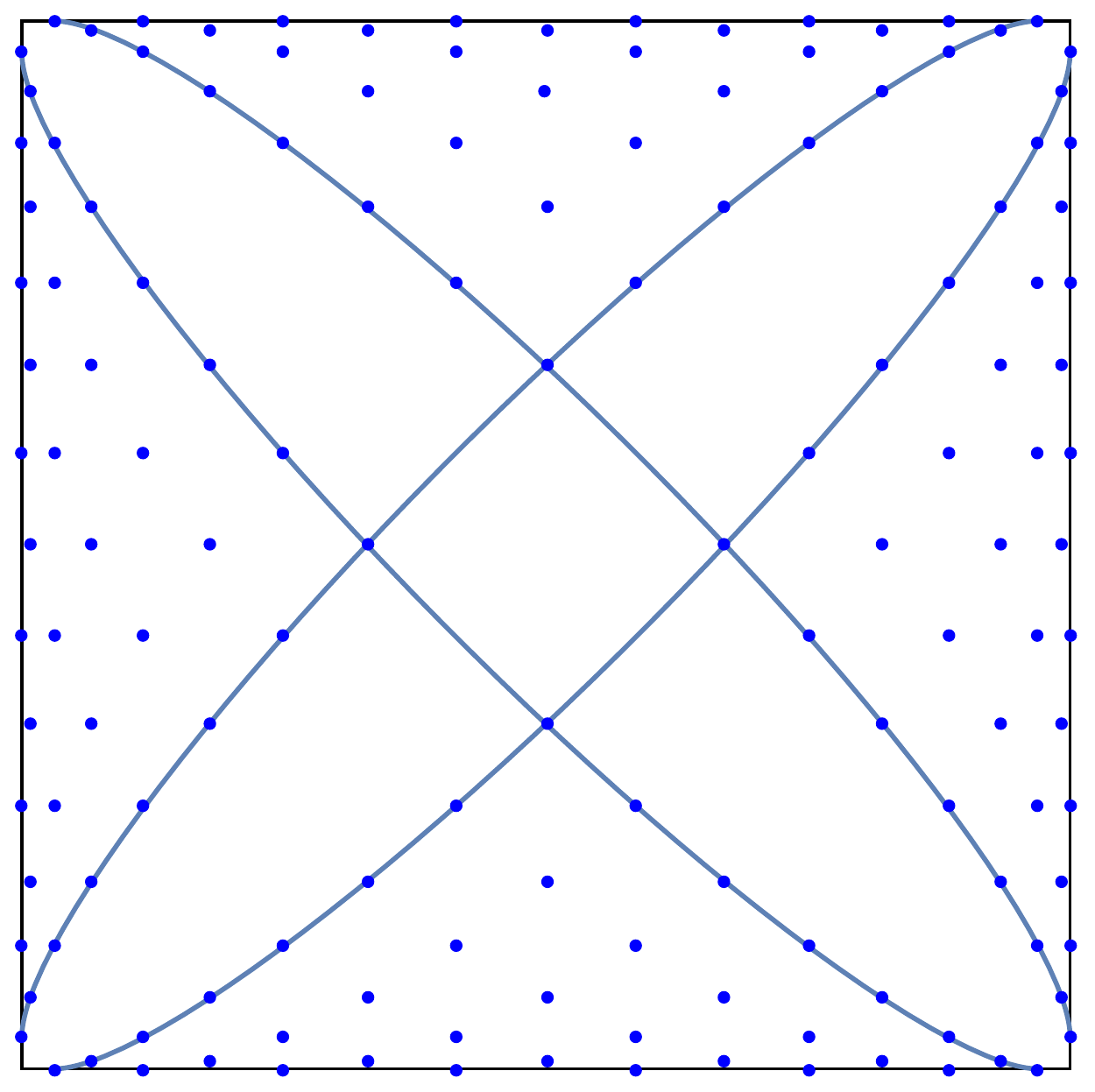} \,\,  \includegraphics[width=2.3in]{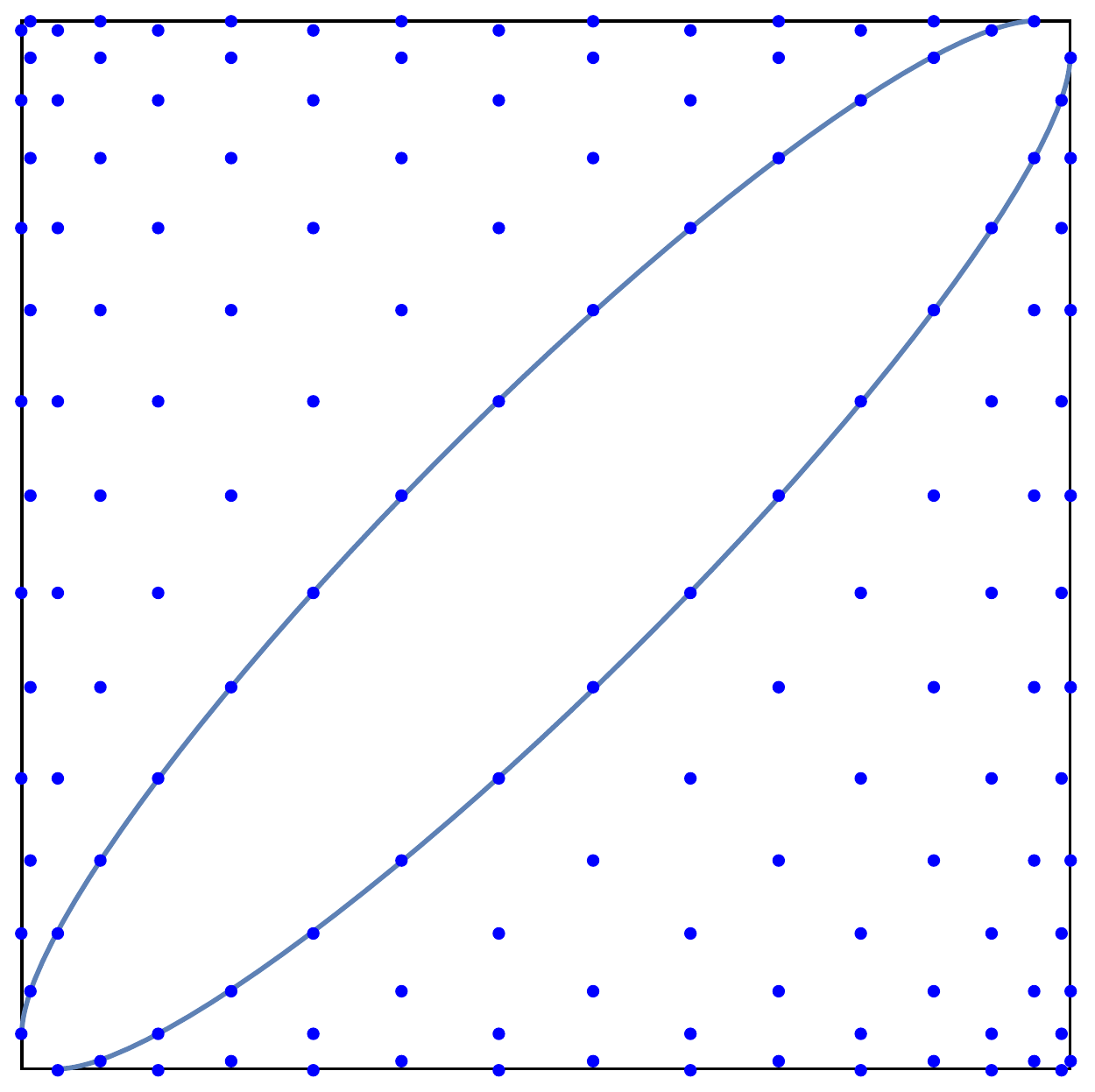} 
   \caption{144 nodes for minimal cubature rule of degree 31 for the weight functions 
   $W_{\frac{1}{2},\frac{1}{2},-\frac{1}{2}}$ (left) and $W_{\frac{1}{2},-\frac{1}{2},-\frac{1}{2}}$ (right)}
 \label{fig:min_curve}
\end{figure}

The weight function $W_{\alpha,\beta, -\frac{1}{2}}$ has a singularity at the diagonal $y=x$ of the square when $\alpha \ne -\frac{1}{2}$, or at
the diagonal $y=-x$ of the square when $\beta \ne -\frac{1}{2}$, or at both diagonals when neither $\alpha$ nor $\beta$ equal to $-\frac{1}{2}$. 
This is reflected in the distribution of the nodes, which are propelled away from these diagonals. Furthermore, for a 
fixed $m$, the points in $X_{2m}$ and $X_{2m+1}$ will be propelled further away for increasing values of $\alpha$ and/or $\beta$.
In Figure \ref{fig:min_curve} we depict the nodes of the minimal cubature rules of degree 31 for $W_{\frac{1}{2},\frac{1}{2},-\frac{1}{2}}$, 
which has singularity on both diagonals, and for $W_{\frac{1}{2},-\frac{1}{2},-\frac{1}{2}}$, which has singularity at the diagonal $y=x$. 
Writing explicitly, these weight functions are 
$$
W_{\frac{1}{2},\frac{1}{2},-\frac{1}{2}}(x,y) = \frac{(x-y)^2(x+y)^2}{\sqrt{1-x^2}\sqrt{1-y^2}} \quad\hbox{and}\quad
 W_{\frac{1}{2}, - \frac{1}{2},-\frac{1}{2}}(x,y) = \frac{(x-y)^2}{\sqrt{1-x^2}\sqrt{1-y^2}}. 
$$
We also depicted the curves that bound the region that does not contain any nodes, which are given in explicit 
parametric formulas in \cite[Proposition 3.6]{X16}. The region without nodes increases in size when $\alpha$ and/or $\beta$ 
increase for a fixed $m$, but they are getting smaller when $m$ increases while $\alpha$ and $\beta$ are fixed. These figures 
can be compared to those in Figure \ref{fig:min} for the case $\alpha= \beta = -\frac{1}{2}$, where the obvious symmetry in $X_n^{\alpha,\beta}$ 
is not evident. 

Let $L_n^{\alpha,\beta} f$ be the interpolation polynomial based on $X_{2m}^{\alpha,\beta}$ when $n=2m$ and on 
$X_{2m+1}^{\alpha,\beta}$
when $n=2m+1$, as defined in \eqref{eq:Lnf}. The asymptotics of the Lebesgue constants for these interpolation polynomials 
can be determined \cite{X12a,X16}.

\begin{thm}
Let $\alpha, \beta \ge -1/2$. The Lebesgue constant of the Lagrange interpolation polynomial $\mathcal{L}_n^{\alpha,\beta}f$ satisfies 
\begin{equation} \label{LebesgueCLn}
  \|\mathcal{L}_n^{\alpha,\beta}\|_\infty = \mathcal{O}(1) \begin{cases} n^{2 \max \{\alpha,\beta\} + 1}, 
  & \max\{\alpha,\beta\} > -1/2, \\
                                (\log n)^2, & \max\{\alpha,\beta\} = -1/2.  \end{cases}
\end{equation}
\end{thm}

It should be mentioned that an explicit formula for the kernel $K_n^*$ in \eqref{eq:Kn*} is known, so that the interpolation
polynomials $L_n^{\alpha,\beta} f$ can be written down in closed form without solving a large linear system of equations.

\section{Minimal Cubature Rules for Constant Weight}

The weight functions in the previous two sections contain the Chebyshev weight functions but do not include the 
weight functions $(1-x^2)^\lambda (1-y^2)^\lambda$ for $\lambda \ne \pm \frac{1}{2}$. In particular, it does not include the 
constant weight  function $W(x,y) =1$. 

For these weight functions, it is possible to establish their existence when $n$ is small. In this section we discuss how
these formulas can be constructed. For cubature rules of degree $2n-2$, we consider the Gaussian cubature rules described 
in the item 2 of Theorem \ref{thm:Gaussian}. For cubature rules of degree $2n-1$, we consider minimal cubature rules 
that attain the lower bound \eqref{eq:3rd-lwbd}. Both these cases can be characterized by non-linear system of equations, which 
may or may not have solutions. We shall describe these equations and solve them for the constant weight function for
small $n$. The known cases for these cubature rules are listed in \cite{Cools1, Cools2}.  

Throughout the rest of this section, we shall assume that $W(x,y) =1$. Let $\mathcal{V}_n$ be the space of orthogonal polynomials 
of degree $n$ with respect to the inner product $\langle f,g\rangle = \frac{1}{2}\int_{\Box} f(x,y) g(x,y) dx dy$. Then an orthonormal 
basis of $\mathcal{V}_n$ is given by
$$
   P_k^n(x,y) = \widehat P_{n-k}(x) \widehat P_k(y), \qquad 0 \le k \le n, 
$$
where $\widehat P_n = \sqrt{2n+1} P_n$ and $P_n$ is the classical Legendre polynomial of degree $n$. In this case, the 
coefficients $B_{n,i}$ in the three-term relations \eqref{eq:3-term} are zero and the three-term relations take the form
\begin{align*}
  x \mathbb{P}_n(x,y) & = A_{n,1} \mathbb{P}_{n+1}(x,y) + A_{n-1,1}^{\mathsf{t}} \mathbb{P}_{n-1}(x,y),  \\
  y \mathbb{P}_n(x,y) & = A_{n,2} \mathbb{P}_{n+1}(x,y) + A_{n-1,2}^{\mathsf{t}} \mathbb{P}_{n-1}(x,y),  
\end{align*}
where $\mathbb{P}_n = (P_0^n, \ldots, P_n^n)^\mathsf{t}$, $A_{n,1}$ and $A_{n,2}$ are given by 
$$
A_{n,1}=\left[ \begin{matrix} a_n&&\bigcirc&0\cr &\ddots&&\vdots \cr
  \bigcirc&&a_0&0\end{matrix} \right]
  \quad \hbox{and} \quad
 A_{n,2} =\left[ \begin{matrix} 0&a_0&&\bigcirc\cr
       \vdots&&\ddots& \cr  0&\bigcirc&&a_n \end{matrix} \right], 
$$
in which 
$$
    a_k:=\frac{k+1}{\sqrt{(2k+1)(2k+3)}}, \qquad k =0,1,2, ... .  
$$

\subsection{Minimal cubature rules of degree $2n-2$}
 
By Theorem  \ref{thm:Gaussian}, the nodes of a Gaussian cubature rule of degree $2n-2$, if it exists, are common zeros
of $\mathbb{P}_n + \Gamma_n \mathbb{P}_{n-1}$ for some matrix $\Gamma_n$ of size $(n+1)\times n$. The latter is characterized in the 
following theorem \cite{X94a}. 

\begin{thm} \label{thm:evenMatrix}
The polynomials in $ \mathbb{P}_n + \Gamma_n \mathbb{P}_{n-1}$ have $n(n+1)/2$ real, distinct zeros if, and only if, $\Gamma_n$ satisfies 
\begin{align} 
  A_{n-1,1} \Gamma_n  = \Gamma_n^\mathsf{t} A_{n-1,1}^\mathsf{t}, \quad & \quad  A_{n-1,2} \Gamma_n
       = \Gamma_n^\mathsf{t} A_{n-1,2}^\mathsf{t}, \label{eq:Ga-1}\\
 \Gamma_n^\mathsf{t} (A_{n-1,1}^\mathsf{t} A_{n-1,2}  - A_{n-1,2}^\mathsf{t} A_{n-1,1}) \Gamma_n  
  &  = (A_{n-1,1} A_{n-1,2}^\mathsf{t} - A_{n-1,2} A_{n-1,1}^\mathsf{t}). \label{eq:Ga-2}
\end{align}
\end{thm}

The equations in \eqref{eq:Ga-1} imply that $\Gamma_n$ can be written in terms of a Hankel matrix $H_n = (h_{i+j})$ of size 
$(n+1) \times n$,
\begin{equation} \label{eq:G-H}
  \Gamma = G_n H_n G_{n-1}^\mathsf{t}, \quad \hbox{where} \quad G_n = \mathrm{diag} \{g_{n,0}, g_{n-1,1},\ldots, g_{1,n-1},g_{0,n}\}
\end{equation}
with 
$$
  g_{n-k,k} = \gamma_{n-k} \gamma_k \quad \hbox{and} \quad \gamma_k = \frac{(2k)! \sqrt{2k+1}}{2^k k!^2}.
$$
Thus, solving the system of equations in Theorem \ref{thm:evenMatrix} is equivalent to solving \eqref{eq:Ga-2} for the
Hankel matrix $H_n$, which is a nonlinear system of equations and its solution may not exist. Since the matrices in both sides of \eqref{eq:Ga-2} are skew symmetric, the nonlinear system consists of $n(n-1)/2$ equations and $2n$ variables. The number of 
variables is equal to the number of equations when $n  = 5$. 

We found the solution when $n= 3, 4, 5$, which gives Gaussian cubature rules of degree $4, 6, 8$. These cases are known in
the literature, see the list in \cite{Cools1}. In the case $n = 3$ and $n = 4$, we were able to solve the system analytically instead
of numerically. For $n = 3$, the matrix $H_3$ takes the form
$$
   H_3 = \frac{4}{27 \sqrt{7}} \left[ \begin{matrix} -\frac{11}{25} & 0 & 1 \\ 0 & 1 & 0 \\ 1 & 0 & \frac{2}{5} \\ 0 & \frac{2}{5} & 0 \end{matrix}\right].
$$
The case $H_4$ is too cumbersome to write down. In the case $n=5$, the system is solved numerically, which has multiple
solutions but essentially one up to symmetry. This solution, however, has one common zero (or node of the Gaussian cubature
rule of degree 8) that lies outside of the square, which agrees with the list in \cite{Cools2}. 

Solving the system for $n > 5$ {\it numerically} yields no solution. It is tempting to proclaim that the Gaussian cubature rules of 
degree $2n-2$ for the constant weight function on the square do not exist for $n \ge 6$, but a proof is still needed. 

\subsection{Minimal cubature rules of degree $2n-1$} Here we consider minimal cubature rules of degree $2n-1$ that attain 
the lower bound \eqref{eq:3rd-lwbd}. By Theorem \ref{thm:op-zeros}, the nodes of such a cubature rule are common zeros of 
$(n+1) - \lfloor n/2\rfloor$ many orthogonal polynomials of degree $n$, which can be written as the elements of 
$U^\mathsf{t} \mathbb{P}_n$, where $U$ is a matrix of size $(n+1) \times (n+1 - \lfloor n/2\rfloor)$ and $U$ has full rank.  

\begin{thm} \label{thm:oddMatrix}
There exist $(n+1) - \lfloor n/2\rfloor$ many orthogonal polynomials of degree $n$, written as $U^\mathsf{t} \mathbb{P}_n$, that have 
$n (n+1) + \lfloor \frac{n}2\rfloor$ real, distinct common zeros if, and only if, $U$ satisfies $U^\mathsf{t} V =0$ for a matrix $V$ of size $(n+1) 
\times \lfloor \frac{n}2\rfloor$ that satisfies 
\begin{align} 
  &  A_{n-1,1} (VV^\mathsf{t}- I) A_{n-1,2}^\mathsf{t}   = A_{n-1,2}(VV^\mathsf{t}- I) A_{n-1,2}^\mathsf{t}, \label{eq:VV-1}\\
   & VV^\mathsf{t} (A_{n-1,1}^\mathsf{t} A_{n-1,2}   - A_{n-1,2}^\mathsf{t} A_{n-1,1}) VV^\mathsf{t} =0, \label{eq:VV-2}
\end{align}
where $I$ denotes the identity matrix. 
\end{thm}
The equation \eqref{eq:VV-1} implies that the matrix $VV^\mathsf{t}$ can be written in terms of a Hankel matrix $H_n$ of size 
$(n+1)\times(n+1)$, 
$$
      VV^\mathsf{t} = I + G_n H_n G_n^\mathsf{t} := W,
$$
where $G_n$ is defined as in \eqref{eq:G-H}. Thus, to find the matrix $V$ we need to solve \eqref{eq:VV-2} for $H_n$ and
make sure that the matrix $W$ is nonnegative definite and has rank $\lfloor \frac{n}2 \rfloor$, so that it can be factored as 
$VV^\mathsf{t}$. The non-linear system \eqref{eq:VV-2} consists of $n(n+1)/2$ equations and has $2n+1$ variables, which  
may not have a solution. 

Comparing with the Gaussian cubature rules of even degree in the previous subsection, however, the situation here is 
more complicated. We not only need to solve \eqref{eq:VV-2}, similar to solving \eqref{eq:Ga-2}, for $H_n$, we also have to 
make sure that the resulting $W$ is non-negative definite and has rank $\lfloor \frac{n}2\rfloor$, which poses an additional 
constraint that is not so easy to verify. 

We found the solutions when $n= 3, 4, 5$ and $6$, which gives minimal cubature rules of degree $5, 7, 9, 11$. These 
cases are all known in the literature, see the list in \cite{Cools1} and the references therein. In the case of $n =4$, there 
are multiple solutions; for example, one solution has all 12 points inside the square and another one has 2 points outside. 
In the case $n = 3, 4, 5$, we were able to solve the system analytically instead of numerically. We give Hankel matrices $H_n$ 
for those cases that have all nodes of the minimal cubature rules inside the square:  
$$
 H_3 = \frac{4}{135} \left[ \begin{matrix} -\frac{8}{35} & 0 & 1 & 0 \\ 0 & 1 & 0 &0 \\ 1 & 0 & 0 & 0 \\  0 & 0 & 0 & \frac{4}{35} \end{matrix}\right],
 \qquad H_4 = \frac{44}{14385} \left[ \begin{matrix} \frac{94}{231} & 1 & 1 & 1 & -\frac{82}{55} \\1 & 1 & 1 & -\frac{82}{55} &1 
     \\ 1 & 1 & -\frac{82}{55} & 1 & 1 \\  1& -\frac{82}{55} & 1 & 1 & 1 \\ -\frac{82}{55} & 1 & 1& 1 &  \frac{94}{231} \end{matrix}\right],
$$
and
$$
 H_5 = \frac{96}{77875} \left[ \begin{matrix} \frac{1151}{2079} & \frac{10 \sqrt{86}}{189} & -\frac{31}{86} & -\frac{1}{9}
  \sqrt{\frac{43}{2}} & 1 & 0 \\
    \frac{10 \sqrt{86}}{189}  & -\frac{31}{86} & -\frac19 \sqrt{\frac{43}{2}} & 1 & 0 &1  \\ 
     -\frac{31}{86} & -\frac19 \sqrt{\frac{43}{2}} & 1 & 0 &1  &\frac19 \sqrt{\frac{43}{2}}\\ 
      -\frac19 \sqrt{\frac{43}{2}} & 1 & 0 &1  &\frac19 \sqrt{\frac{43}{2}} & -\frac{31}{86} \\
       \\ 1 & 0 &1  &\frac19 \sqrt{\frac{43}{2}} & -\frac{31}{86} & -  \frac{10 \sqrt{86}}{189}  \\
       0 &1  &\frac19 \sqrt{\frac{43}{2}} & -\frac{31}{86} & -  \frac{10 \sqrt{86}}{189} & \frac{1151}{2079}
       \end{matrix}\right].
$$
Once $H_n$ is found, it is easy to verify that $W$ satisfies the desired rank condition and is non-negative definite. We can 
then find $U$, or the set of orthogonal polynomials, and then find common zeros. For example, when $n =5$, we have 
$4$ orthogonal polynomials of degree $5$ given by 
\begin{align*}
 Q_1(x,y)= & \frac{10 \sqrt{86}}{189}  P_0^5(x,y) +   \frac{1081 \sqrt{11}}{2835 \sqrt{3}}  P_1^5(x,y) +   P_5^5(x,y),\\
 Q_2(x,y)= & \frac{205}{21 \sqrt{33}}  P_0^5(x,y) +  \frac{10 \sqrt{86}}{189}  P_1^5(x,y) +  P_4^5(x,y), \\
 Q_3(x,y)= & -  \frac{5 \sqrt{438}}{27 \sqrt{77}} P_0^5(x,y) +  \frac{62 \sqrt{5}}{81\sqrt{21}}  P_1^5(x,y) +  P_3^5(x,y), \\
 Q_4(x,y)= & -  \frac{10 \sqrt{5}}{3 \sqrt{77}} P_0^5(x,y) -  \frac{\sqrt{430}}{9\sqrt{21}}   P_1^5(x,y) +   P_2^5(x,y),
\end{align*}
which has 17 real common zeros inside the square. Only numerical results are known for the case $n=6$. We also 
tried the case $n =7$, but found no solution numerically. 

\medskip\noindent
{\bf Acknowledgements}\quad The author thanks two anonymous referees for their careful reading and corrections.

\end{document}